\documentclass[12pt]{amsart}
\usepackage{amsfonts, amssymb,amsmath, amsthm, amsxtra, latexsym, amscd}
\usepackage[latin1]{inputenc} 

\def\draft{\centerline{(Draft {\the \day}/{\the\month} \the \year.)}}

\theoremstyle{definition}
\newtheorem{theo+}    {Theorem}      [section]
\newtheorem{prop+}  [theo+]  {Proposition}
\newtheorem{coro+}  [theo+]  {Corollary}
\newtheorem{lemm+}  [theo+]  {Lemma}
\newtheorem{deep+}  [theo+]  {Deep Result}
\newtheorem{fact+}  [theo+]  {Fact}
\theoremstyle{definition}
\newtheorem{exam+}  [theo+]  {Example}
\newtheorem{rema+}  [theo+]  {Remark}
\newtheorem{defi+}  [theo+]  {Definition}
\newtheorem{xca+}[theo+]{Exercise}

\numberwithin{equation}{section}

\def\draft{\centerline{(Draft {\the \day}/{\the\month} \the \year.)}}
\bigskip

\def\refn#1.#2{\expandafter\def\csname#1\endcsname{[#2]}}
\def\refnr#1.{\csname#1\endcsname}

\def\fa{\mathfrak a}

\def\fg{\mathfrak g}
\def\fk{\mathfrak k}
\def\fh{\mathfrak h}
\def\fl{\mathfrak l}

\def\fn{\mathfrak n}
\def\fp{\mathfrak p}
\def\fq{\mathfrak q}

\def\a{\alpha}

\def\Claminv2{|C(\Lambda)|^{-2}}
\def\ga{\gamma}

\def\varepsi{\varepsilon}
\def\lam{\lambda}

\def\Ome{\Omega}

\def\de{d\varepsilon}

\def\Aa2D{A^{\a,2}(D)}
\def\bAa2D{\overline{A^{\a,2}(D)}}
\def\Ab2D{A^{\beta,2}(D)}
\def\bAb2D{\overline{A^{\beta,2}(D)}}

\def\Norm#1_#2{\Vert#1\Vert_{#2}}

\def\phipl12{\phi_{p_{l_1}, p_{l_2}}}
\def\phip01{\phi_{p_{0}, p_{0}}}

\def\a{\alpha}

\def\ga{\gamma}
\def\Claminv2{|C(\Lambda)|^{-2}}

\def\varepsi{\varepsilon}

\def\sig{\sigma}
\def\Sig{\Sigma}

\def\lam{\lambda}

\def\exp{\operatorname{exp}}

\def\span{\operatorname{span}}

\def\de{d\varepsilon}

\def\Aa2D{A^{\a,2}(D)}
\def\bAa2D{\overline{A^{\a,2}(D)}}
\def\Ab2D{A^{\beta,2}(D)}
\def\bAb2D{\overline{A^{\beta,2}(D)}}

\def\phipl12{\phi_{p_{l_1}, p_{l_2}}}
\def\phip01{\phi_{p_{0}, p_{0}}}
\def\m{\underline{\bold m}}
\def\n{\underline{\bold n}}

\def\bc{\mathbb C}

\def\alg/{algebra} 
\def\Alg/{Algebra} 
\def\alt/{alternative} 
\def\anal/{analytic}
\def\analfunc/{\anal/\ \func/}
\def\Ans/{\it Answer. \normal}
\def\ass/{associative}
\def\nass/{non-\ass/}
\def\autom/{automorphism}
\def\homom/{homomorphism}
\def\isom/{isomorphism}
\def\bdd/{bounded}
\def\Bdd/{Bounded}
\def\bddsymdom/{bounded \sym/ \dom/}
\def\Cartdom/{Cartan \dom/}
\def\bdry/{boundary}
\def\bsd/{\bdd/ \symdom/}
\def\bv/{boundary value}
\def\cf/{{\it cf}\.}
\def\Cf/{{\it Cf}\.}
\def\charr/{character}
\def\coeff/{coefficient}
\def\comm/{commutative}
\def\cpct/{compact}
\def\compl/{complex}
\def\comp/{complex}
\def\Comp/{Complex}
\def\conf/{conformal}
\def\conj/{conjugate}
\def\conn/{connect}
\def\cont/{continuous}
\def\conv/{converge} 
\def\convc/{convergence}
\def\convt/{convergent}
\def\convx/{convex}
\def\coord/{coordinate}
\def\lcoord/{local coordinate}
\def\Corr/{Corresponding}
\def\corr/{corresponding}
\def\corrd/{correspond}
\def\cov/{covariant}
\def\decomp/{decomposition}
\def\deco/{decompose}
\def\diff/{different} 
\def\Diff/{Different} 
\def\dimn/{dimension} 
\def\distr/{distribution} 
\def\div/{diverge} 
\def\dom/{domain}
\def\eg/{\hbox{\it e.g}\.}
\def\eigenf/{eigen\-\func/}
\def\eigensp/{eigen\-space}
\def\eigenv/{eigen\-value}
\def\eq/{equation}
\def\equa/{equation}
\def\de/{\diff/ial \equa/}
\def\do/{\diff/ial operator}
\def\ode/{ordinary \de/}
\def\pde/{partial \de/}
\def\pdo/{partial \diff/ial operator}
\def\psdo/{pseudo \diff/ial operator}
\def\fin/{finite}
\def\Ex/{\it Example.\ \normal}
\def\Exnr#1/{\it Example #1.\ \normal}
\def\foll/{follow}
\def\follg/{following}
\def\Follg/{Following}
\def\func/{function}
\def\Func/{Function}
\def\Fonc/{Fonc\-tion}
\def\fonc/{fonc\-tion}
\def\Funk/{Funk\-tion}
\def\funk/{Funk\-tion}
\def\gen/{general}
\def\har/{harmonic}
\def\Hint/{\it Hint. \normal}
\def\hist/{historic}
\def\histcl/{historical}
\def\hol/{holo\-morphic}
\def\homog/{ho\-mo\-ge\-ne\-ous}
\def\hyp/{hyper\-bolic}
\def\hyperg/{hyper\-geometric}
\def\ie/{\hbox{\it i.e.}}
\def\iff/{if and only if}
\def\ineq/{inequality}
\def\infra/{{\it inf\-ra}}
\def\ultra/{{\it ult\-ra}}
\def\Inpart/{In particular}
\def\inpart/{in particular}
\def\instof/{instead of}
\def\interps/{interpolation space}
\def\interp/{interpolation}
\def\Interp/{Interpolation}
\def\interpr/{Interpretation}
\def\Intr/{Introduction}
\def\intv/{interval}
\def\inv/{invariant}
\def\invc/{invariance}
\def\Iowords/{In other words}
\def\iowords/{in other words}
\def\ipr/{inner product}
\def\irred/{irreducible}
\def\lb/{line bundle}
\def\lin/{linear}
\def\lhs/{left hand side}
\def\rhs/{right hand side}
\def\loc/{local}
\def\math/{mathematic} 
\def\mathcn/{\math/ian}
\def\manif/{manifold}
\def\meas/{measure}
\def\measl/{measurable}
\def\mero/{mero\-morphic}
\def\mon/{monomial}
\def\monog/{monogenic}
\def\mult/{multiple}
\def\multy/{multiply}
\def\multn/{multiplication}
\def\nas/{necessary and sufficient}
\def\nbd/{neighborhood}
\def\neg/{negative}
\def\nondeg/{nondegenerate}
\def\Oohand/{On the other hand}
\def\oohand/{on the other hand}
\def\Oonhand/{On the one hand}
\def\oonhand/{on the one hand}
\def\oper/{operator}
\def\orth/{ortho\-gonal}
\def\orthon/{ortho\-normal}
\def\otoh/{on the other hand}
\def\quat/{quaternion}
\def\pp/{\hbox{a. e.}}
\def\psh/{plurisubharmonic}
\def\pol/{polynomial}
\def\pot/{potential}
\def\pos/{positive}
\def\princ/{principle}
\def\prob/{probability}
\def\proj/{projective}
\def\projn/{projection}
\def\Proof/{\it Proof:\normal}
\def\Rem/{\it Remark\normal}
\def\Remnr#1/{\it Remark\ \normal #1. }
\def\rep/{representation}
\def\meta/{metaplectic representation}
\def\repr/{reproducing}
\def\reprker/{reproducing kernel}
\def\resp/{respective} 
\def\resply/{respectively}
\def\restr/{restriction}
\def\sa/{self-adjoint}
\def\st/{such that}

\def\sol/{solution}
\def\ru/{space}
\def\sph/{spherical}
\def\ssp/{sub\ru/}
\def\sym/{symmetric}
\def\Sym/{Symmetric}
\def\symb/{symbol}
\def\symbc/{symbolic}
\def\symdom/{\sym/ domain}
\def\symp/{symplectic}
\def\Theor#1/{\fet Theorem #1.\ \normal}
\def\Lem#1/{\fet Lemma #1.\ \normal}
\def\Lemma/{\fet Lemma.\ \normal}
\def\topl/{topology}
\def\topll/{topological}
\def\transf/{transform}
\def\transl/{translation}
\def\transfn/{transformation}
\def\transv/{transvectant}
\def\trig/{trigonometric}
\def\tril/{trilinear}
\def\trilf/{trilinear form}
\def\uhp/{upper halfplane}
\def\uhs/{upper halfspace}
\def\vb/{vector bundle}
\def\vf/{vector field}
\def\vsp/{vector space}
\def\wrt/{with respect to}
\def\Wlog/{Without loss of generality}
\def\a{\alpha}

\def\lam{\lambda}

\def\sig{\sigma}
\def\Sig{\Sigma}

\def\Ab/{Abel}
\def\Ban/{Banach}
\def\Bansp/{\Ban/ space}
\def\Belt/{Bel\-tra\-mi}
\def\Berg/{Berg\-man}
\def\Bern/{Ber\-nou\-lli}
\def\Berz/{Berezin}
\def\Bess/{Bessel}
\def\Cart/{Car\-tan}
\def\Cay/{Cay\-ley}
\def\CG/{Clebsch-Gordan}
\def\Cl/{Clifford}
\def\CR/{Cauchy-Rie\-mann}
\def\Dir/{Dirichlet}
\def\Eucl/{Euclide}
\def\F/{Fourier}
\def\Hank/{Hankel}
\def\Hankf/{\Hank/ form}
\def\Herm/{Hermite}
\def\Hilb/{Hilbert}
\def\Hilbs/{Hilbert space}
\def\Hilbsp/{Hilbert space}
\def\HS/{Hilbert-Schmidt}
\def\Lag/{La\-grange}
\def\Lap/{La\-place}
\def\LapBelt/{\Lap/-\Belt/}
\def\Leb/{Lebesgue}
\def\Marc/{Mar\-cin\-kie\-wicz}
\def\Moeb/{Moebius}
\def\Moebt/{Moebius transformation}
\def\Moebtransfn/{Moebius transformation}
\def\Pla/{Plan\-che\-rel}
\def\Poin/{Poin\-car\'e}
\def\Riem/{Rie\-mann}
\def\Riemn/{\Riem/ian}
\def\psRiemn/{pseudo-\Riem/ian}
\def\Riems/{Rie\-mann surface}
\def\Schroe/{Schr\"odinger}
\def\Weier/{Weier\-strass}
\def\anal/{analytic}
\def\bsd/{bounded symmetric domain  }
\def\bdd/{bounded}
\def\calc/{calculation}\def\conj{conjugate}
\def\calci/{calculating}\def\eg{e.g.}
\def\conj/{conjugate}
\def\deco/{decomposition}
\def\eg/{e.g.}
\def\fct/{function}
\def\gp/{group}
\def\hw/{highest weight}
\def\hwv/{highest weight vector}
\def\hwvs/{highest weight vectors}
\def\lw/{lowest weight}
\def\lwv/{lowest weight vector}
\def\lwvs/{lowest weight vectors}
\def\hds/{holomorphic discrete series}
\def\iff/{if and only if}
\def\inv/{invariant}
\def\irrde/{irreducible decomposition}
\def\meas/{measure}
\def\transf/{transform}
\def\rep/{representation}
\def\resp/{respectively}
\def\inters/{intertwines}
\def\interg/{intertwining}
\def\meta/{metaplectic representation}
\def\qu/{quaternion}
\def\rep/{representation}
\def\symdom/{ symmetric domain}
\def\st/{such that}
\def\shd/{subhead}
\def\transf/{transform}
\def\wrt/{with respect to}

\def\Norm#1#2#3{\Vert#1\Vert^{#3}_{{#2}Â¥}}

\def\Res{\operatorname{Res}}

\begin{document}
\title[Principal series representations and holomorphic extensions]
{Degenerate principal series representations
and their holomorphic extensions}
\author{Genkai Zhang}
\address{Mathematical Sciences, Chalmers University of Technology and
Mathematical Sciences, G\"oteborg University
 , SE-412 96 G\"oteborg, Sweden}
\email{genkai@math.chalmers.se}
\thanks{Research  partially supported by the Swedish
Science Council (VR), the Max-Planck institute of mathematics, Bonn Germany}

\def\bbK{\mathbb K}
\def\bKk{\mathbb K^k}
\def\bKn{\mathbb K^n}
\def\Cos^2{\text{Cos}^2}

\begin{abstract} Let $X=H/L$ be an irreducible real bounded symmetric
domain  realized as a real form
in an Hermitian
symmetric domain
$D=G/K$. The intersection $S$ of the
 Shilov boundary of $D$ with $X$ defines
a distinguished subset  of the topological boundary  of $X$ and is invariant
under $H$ and can also be realized as $S=H/P$ for certain parabolic subgroup $P$ of
$H$. We  study the spherical
representations $Ind_P^H(\lam)$ of $H$ induced from $P$.
We find  formulas for the spherical functions
 in terms of the Macdonald ${}_2F_1$ hypergeometric
function. This
generalizes the earlier result of  Faraut-Koranyi for
Hermitian symmetric spaces $D$.
We consider a class of $H$-invariant integral intertwining
operators from the representations
$Ind_P^H(\lam)$ on  $L^2(S)$ to the holomorphic
representations of $G$ on $D$ restricted to $H$.
We construct a new class of complementary
series for the groups $H=SO(n, m)$,  $SU(n, m)$ 
(with $n-m >2$) and  $Sp(n, m)$ (with $n-m>1$).  We realize them
as a discrete component in the branching rule of
the analytic continuation of the
holomorphic discrete series of $G=SU(n, m)$,
$SU(n, m)\times SU(n, m)$
and $SU(2n, 2m)$ respectively. 
\end{abstract}
\maketitle









\def\Dc{\mathbb D}
\def\Vc{V_{\mathbb C}}
\def\Bnu{B_{\nu}} 
\def\el{e_{\blam}} 
\def\bnul{b_{\nu}(\lam)} 
\def\GGa{\Gamma_{\Ome}} 
\def\fc#1#2{\frac{#1}{#2}} 
\def\SS{\mathcal S}

\def\mP{\mathcal P}
\def\mF{\mathcal F}
\def\SdN{\mathcal S^N}
\def\eSdN{\mathcal S_N^\prime}
\def\Gc{G_{\mathbb C}}
\def\rqs{G\backslash K_{\mathbb C}}

\baselineskip 1.25pc

\section{Introduction}

\def\Pn{\mathcal P_{\n}}
\def\Hnu{\mathcal H_{\nu}}
\def\HtH{\mathcal H_{\nu}  \otimes
\overline{\mathcal H_{\nu}}}
  
Since the work of
Kashiwara and Vergne \cite{KV}
on the tensor product decomposition of metaplectic representations
and of  Howe \cite{Howe-jams} on the dual pair correspondence \cite{Howe-reciprocity},
there has been  intensive study  on the 
branching rule of minimal and singular
 representations under various subgroups;
see e.~g.~\cite{Kobayashi-icm02, KO-adm}
and references therein.
The purpose of the present paper is to 
find certain irreducible discrete parts
 of the branching 
of  scalar 
holomorphic representations $\pi_{\nu}$  
of an Hermitian Lie group $G$ of higher rank 
under a
symmetric subgroup $H$, with the Riemannian 
symmetric domain $X=H/L$ being
a real form of the Hermitian symmetric domain $D=G/K$.
For larger parameter $\nu$ it is equivalent \cite{gkz-bere-rbsd}
to the regular action of $H$ on $L^2(X)$, whose decomposition is well-known \cite{He2}
and is a continuous sum of the principal series
representations induced on the minimal Iwasawa parabolic subgroup.
However for smaller parameter $\nu$ the decomposition
is rather complicated with the  continuous parts
being integration over various hyperplanes in
the complex dual of the real Cartan subalgebra and with discrete parts,
and the full decomposition is not known; see \cite{Neretin-plan-beredef}
for some examples. It is thus worthwhile to find the discrete
components which in certain sense are the most interesting part.
A pivot example of such cases is when $X$ is
itself a complex bounded symmetric domain $G/K$ realized
as a diagonal part in the domain $X\times \bar X=G\times G/K\times K$.
In this case the branching rule above is
the tensor product decomposition 
of $\pi_{\nu}\otimes \overline{\pi_{\nu}}$
of the holomorphic representation with its complex conjugate.
However for smaller parameter $\nu$ there might have
some discrete components. A full decomposition
is done in \cite{gz-tensormin} for $\nu$ being the last
Wallach point, and in this case
there are finitely many discrete
components appearing
for the  non-tube type domain 
$SU(r, r+b)/S(U(r)\times U(r+b))$,
$r>1$, $b>2$, and they  are then some complementary
series representations. We thus get a realization
of them in the space $\HtH$, namely the space of
Hilbert-Schmidt operators on  $\mathcal H_{\nu}$, which can be viewed
as a quantization of the complementary series representation; this
is another part of our motivation.
In the  rank one case
the appearance of the discrete
parts 
in the tensor product
 $\pi_{\nu}\otimes \overline{\pi_{\nu}}$ for the  group
$SU(n, 1)$ has been studied 
in \cite{Vershik-G-G}
for $n=1$ and 
 for $n>1$ in  \cite{ehprz-comser};
similar results hold for the branching of
holomorphic representation of $SU(n, 1)$ 
under $SO(n, 1)$ 
and 
$Sp(\frac n2, 1)$ \cite{Dijk-Hille-jfa}.
Certain examples of higher rank cases have also been
studied earlier in \cite{oz-cjm} and
\cite{Neretin-plan-beredef}. 
In \cite{Neretin-Olshanski} 
Neretin and Olshanski  discovered
that certain complementary series
representation for $SO(p, q)$ appears in the branching
of the minimal representation
of $SU(p, q)$; see also \cite{HS-so}.
Here we give a systematic study
of the appearance of complementary series
in the branching of the holomorphic representation
$\pi_{\nu}$ for all groups $G$ and $H$.

We will consider 
representations of $H$ 
that are induced from some maximal parabolic subgroup $P$.
 More precisely  we consider the boundary
$S=H/P$  defined as the  intersection 
of the Shilov boundary of $D$ with the topological
closure of $X$, which we may call the real Shilov boundary
of $X$.
The corresponding spherical functions on $X=H/L$
can be realized as certain Poisson integral, and they
have natural analytic continuation
to the whole domain $D$. We find
first an expansion formula
for the spherical functions in terms
of the $L$-invariant polynomials on $D$, which
are certain hypergeometric functions
studied earlier by Faraut-Koranyi and Macdonald.
For that purpose we generalize our earlier
results \cite{gz-br2}
about characterizing the $L$-invariant polynomials
in terms of the Jack symmetric  polynomials
 to non-tube type domains,
and we find their Fock space norm using
the result of Dunkl \cite{Dunkl-cmp-98};
combining this with the result of Faraut-Koranyi 
\cite{FK} we
find then their Bergman space norm. Even though
the computations are rather technical
the end results on the spherical functions
are simple and appealing; see Theorem 5.1.

We study consequently the question
of describing those spherical representations
that appear discretely in the branching
of the holomorphic representations
of $G$ under $H$. It turns out this
happens only, roughly speaking,  when $\nu$
is smaller  and when $H=SO(n, m)$,  $SU(n, m)$ 
 and  $Sp(n, m)$ with $n-m$ sufficiently larger
than $\nu$.
The result of this type is heuristically plausible as the 
the spaces of singular
holomorphic representations 
are certain Sobolev type space and may contain $H$-invariant
subspaces holomorphic functions
with  boundary values (called also trace in classical analysis)
on certain boundary components of $X=H/L$, and
the subspaces may form certain complementary series
representation; conversely  complementary
series representations of $H$ might usually
be realized  on space of functions
on boundary components of $X$ 
with extra smoother property (see e.g.
\cite{Kirillov-bams-99} for some precise
statements and related conjectures)
and they may have holomorphic continuation 
on $D$. However
the precise statement is rather subtly.
The parameter for the induced representations are outside
the unitary range so that
we may call them complementary series;
in rank one
case they are precisely the known
complementary series. 
We prove that when $\nu$ in the continuous part
of the Wallach set and is small compared with the root
multiplicity $b$ there appear discrete parts in
the branching of the holomorphic representation
of $G$ under $H$.

Our results on realization of the complementary
series representation as discrete components
in the holomorphic representations
can also be viewed as  representation
theoretic study of the holomorphic
extension of spherical functions. Roughly we are mostly 
concerned with constructing unitary spherical 
representation of $H$ so that its holomorphic extension
is a discrete component in the unitary holomorphic
representation of the larger group $G$. For a general
Riemannian symmetric space $H/L$ the
holomorphic extension of spherical
functions is of considerable interests; see e.g. \cite{SK-2}.
Also Kobayashi \cite{Kobayashi-vis-act} has recently
 introduced some geometric concept characterizing
the multiplicity free branching rule.

{\it Acknowledgment.} I would like to thank Toshiyuki
Kobayashi, Siddhartha Sahi and Henrik Sepp\"anen for
some helpful discussions. I would also like to thank the support and the hospitality of the Max-Planck institute a
and Haussdorff institute in Mathematics, Bonn, Germany, during my stay in July 2007 where part of this work was finished.

For the reader's convenience we list the main symbols used in the paper.

\begin{enumerate}

\item  $D =G/K$, a 
bounded symmetric  type domain
in a complex vector space $V_{\bc}=V^{\mathbb C}$
of rank $r^\prime$ with
root multiplicity $(1, a^\prime, 2b')$.

\item  $X =H/L$, an irreducible
real bounded symmetric  type domain
in a real form $V\subset V^{\mathbb C}$
of rank $r$  with root multiplicity $(\iota-1, a, 2b)$.

\item  $S= Le=H/P$ a distinguished boundary component
in the topological boundary of $X$.

\item $\mathfrak{h} =\mathfrak{l }+\mathfrak{q}$, 
the Cartan
decomposition of $\mathfrak g$,
 $\fa\subset \mathfrak h$ a maximal abelian
subspace of $\fp$, $\Sig(\fh, \fa)$ the
root system.

\item  $h(z, \bar w)$ an irreducible polynomial on $\Vc\times\overline{\Vc}$,
the Bergman reproducing kernel is $h(z,\bar w)^{-p}$
with $p=a^\prime(r^\prime -1) +2 +b^\prime$ the genus of domain $\mathbb D$.

\item For  a tuple (partition) $\m=(m_1, \cdots, m_r)$
with $m_j\in \mathbb N$ (non negative  integers) and
$m_1\ge \cdots \ge m_r\ge 0$ the generalized Pochammer symbol
is
$$
(c)_{\m, \beta}=\prod_{j=1}^{r} (c-\beta(j-1))_{m_j}=
\prod_{j=1}^{r} \prod_{k=1}^{m_j}
(c-\beta(j-1) +k-1).
$$

\end{enumerate}

\section{Preliminaries}

\subsection{Bounded symmetric domains  $D$ and holomorphic representations}

We recall very briefly in this and next subsections
some preliminary results on bounded
symmetric domains and fix  notation;
see e.g.  \cite{FK, Loos-bsd} and references therein.

Let $D=G/K$ be as in the previous section
an irreducible bounded symmetric domain
  in a $d$-dimensional
complex vector space $V_{\mathbb C}=\bc^d$
of rank $r^\prime$. (The symbol $r$ will
be reserved for the rank of the real bounded symmetric
domain $X$ in next subsection.)  
Let $\fg=\fk+\fp$
be the Cartan decomposition
and  $\fg^{\bc}=\fp^{+}+\fk^{\bc}+\fp^-$
be the Harish-Chandra decomposition of its
complexification.
Let 
$(1, a^\prime, 2b^\prime)$ be the root multiplicities
(or Peirce invariants in terms of the Jordan triple)
of the real root system of $\fg$, with $1$
being that of the longest roots, $2b^\prime$
 of the shortest roots, and $a$ the middle.
The rank $r^\prime$ and the multiplicities
form a quadruple characterizing $D$
\begin{equation}
  \label{eq:D-ch}
\boxed{D: (r^\prime, 1, a^\prime, 2b^\prime)  }
\end{equation}
and will be compared with that for the subdomains $X$ below.
According to the classification \cite{He1}
of bounded symmetric domains
the possible value of
$(r^\prime,  a^\prime, 2b^\prime) $
is
\begin{equation}
  \label{eq:D-cl}
(r^\prime, 2, 2b), 
\quad (r^\prime,  1, 0), \quad (r^\prime,  4, 0), 
 \quad (r^\prime,  4, 2), 
 \quad (2,  a^\prime, 0),
 \quad (2,  6, 8),  \quad (3,  8, 0).
\end{equation}

The space $V_{\mathbb C}=\fp^+$
has  a Jordan triple structure so that
the subspace $\fp$ of the Lie algebra $\fg$,
when realized as a space of holomorphic vector
fields on $\mathbb D$, consists of
vector fields of the  form
\begin{equation}
\label{p-form}
\xi_v=\xi_v(z)=v-Q(z)\bar v, \quad v\in V_{\bc},
\end{equation}
where $Q(z): \bar {V_{\bc}} \mapsto
 {V_{\bc}}$ is quadratic in $z$.
 We denote $
\{x\bar y z\}=D(x, \bar y)z$ 
the Jordan triple product
$$ 
\{x\bar y z\}=D(x, \bar y)z=
(Q(x+z)- 
Q(x)- Q(z))\bar y.
$$
We fix
a $K$-invariant Hermitian inner product
 $(\cdot, \cdot)$ on $V_{\mathbb C}$ so that a minimal
tripotent has norm $1$. We let $dm(z)$ be
the corresponding Lebesgue measure. The Bergman
reproducing kernel is 
up to a positive constant
of the form $h(z, \bar w)^{-p}$ where $p$ is the genus
of $\mathbb D$, defined by $p=a(r^\prime-1) +2 +b^\prime$,
 and $h(z, \bar w)$ is an irreducible
polynomial holomorphic in $z$ and anti-holomorphic
in $w$. In particular the function $h(z, \bar w)$ satisfies
the following transformation property under the group 
$G$,
\begin{equation}
\label{h-trans}
h(gz, \overline {gw})= 
J_{g}(z)^{\frac 1 p} 
h(z, \bar w)\overline{J_{g}(w)^{\frac 1 p} }, \quad g\in G,
\end{equation}
where $J_g$ is the Jacobian of the holomorphic mapping $g$.

We denote by $\mathcal F(V_{\bc})$ 
 the Fock space
of entire functions on $V_{\bc}$. 
Let $\nu>p-1=a(r^\prime-1) +1 +b$ and  consider the probability measure
$d\mu_{\nu}(z)=c_{\nu}^\prime h(z,\bar z)^{\nu -p}dm(z)$
where $c_{\nu}^\prime$
is the normalization constant,
and the corresponding weighted Bergman space $\Hnu=\Hnu(\mathbb D)$
of holomorphic functions $f$
so that
$$
\Vert f\Vert_{\nu}^2=\int_{\mathbb D}|f(z)|^2 d\mu_{\nu}(z)<\infty.
$$
It has reproducing
kernel $h(z, \bar w)^{-\nu}$.
The group $G$ acts unitarily on
$\Hnu$ via the following
\begin{equation}
\label{pi-nu}
\pi_{\nu}(g) f(z)= J_{g^{-1}}(z)^{\frac \nu p}f(g^{-1}z),  
\end{equation}
and it forms a unitary projective representation of $G$.
We let $\mathcal O(D)$ be the space of all holomorphic
functions on $D$. 
The formula 
(\ref{pi-nu})
defines also a representation of $G$ on the space $\mathcal O(D)$.

The representation has an analytic continuation in $\nu$
and the whole set of $\nu$ so that it still defines
an irreducible unitary representation on a proper subspace
of holomorphic functions is given by the so-called
Wallach set
$$
 W=\{0, \frac {a^\prime} 2, \cdots, \frac {a^\prime} 2 (r^\prime-1)\}
\cup (\frac a2(r^\prime-1), \infty).
$$
The corresponding Hilbert space for $\nu\in  W$
 will also be denoted by $\Hnu$
and the norm by $\Vert \cdot\Vert_\nu$. The discrete
points in the set will be also referred as singular
Wallach point (to differ the discrete component
in the branching rule).

We summarize some related results  (see e.g. \cite{FK}) in the following 
\begin{theo+}  Let  $D=G/K$ be as  above.  The space $  \mathcal P$
of holomorphic polynomials on $V_{\bc}$ decomposes into irreducible
subspaces under $K$, with multiplicity one as: 
\begin{equation}
\label{hua-sch}
  \mathcal P\cong \sum_{\underline{\mathbf{n}} \ge 0} \Pn.
\end{equation}
Each $\Pn$ is of lowest weight $-\underline{\mathbf{n}}
=-(n_1\gamma_1 +\cdots +n_{r^\prime} \gamma_{r^\prime})$
 with $n_1\ge \dots \ge n_{r^\prime}\ge 0$ 
and $\ga_1 >\cdots >\ga_r$ the Harish-Chandra strongly orthogonal
roots.
For each nonzero $f\in \Pn$  it holds
$$
{\Vert f\Vert_{\mathcal F}^2}
={(\nu)_{\n, \frac{a^\prime}2}}
{
\Vert f\Vert_{\nu}^2
},
$$
where $(\nu)_{\n, \frac {a'}2}$
is the generalized Pochammer symbol in Section 1. The reproducing
kernel $h(z, \bar w)^{-\nu}$ has the
following expansion
\begin{equation}
\label{FK-exp}
h(z, \bar w)^{-\nu}=\sum_{\n}(\nu)_{\n, \frac{a^\prime}2}K_{\n}(z, \bar w),
\end{equation}
where $
K_{\n}(z,\bar  w)$ is the reproducing kernel of
$\Pn$ in the Fock space.
In particular for 
$\nu=\frac {a^\prime}2 (j-1)$, $1\le j\le r^\prime$, we have
\begin{equation}
  \label{eq:sig-wa-hil}
h(z, \bar w)^{-\nu}=\sum_{\n; n_j=0} (\nu)_{\n}K_{\n}(z, \bar w), \quad
\mathcal H_{\nu}= \sum_{\n: n_{j}=0}^{\oplus}\mathcal P_{\n}.
\end{equation}
\end{theo+}

\subsection{Real forms $X$ of $D$}

Let $V\subset V_{\bc}$ be a real form of $V_{\bc}$,
$V_{\bc}=V+i V$ and let $X=V\cap D$ be
the corresponding real form of $D$. $X$ is
called a real bounded symmetric domain
if the real involution with respect to $V$ preserves
the domain $D$. In this case $V$ is a real Jordan triple
and $X$ is a Riemannian
symmetric space, $X=H/L$,
with induced metric
from that of $D$, 
realized as a bounded domain in $V$. 
Here we take $H$ the
connected component of the subgroup of $G$ preserving
$X$.
The  most well-studied case is when $H/L$ 
is the symmetric cone  in the Siegel tube domain  $G/K$, namely Type A below.
To have
a some what unified treatment
we exclude  the rank one case and we will only consider those  irreducible  $X$.

 Let $\fh$ be the
Lie algebra of $H$ and 
 $\fh=\fl +\fq$ be the Cartan decomposition.
We let $\fa$ be a maximal abelian subspace of
$\fq$.  We fix a frame of minimal tripotents
 $\{e_1, \cdots, e_r\}$ of $V$.
 The corresponding
vector fields  (see \ref{p-form}),  
$\xi_j=\xi_{e_j}$, $j=1, \cdots, r$ 
form  a maximal abelian subspace  $\fa$
of $\fq$.  We will also view $\fa$ as a subspace of $V$. 
We let
\begin{equation}
  \label{eq:max-tri}
  e:=e_1 +\cdots +e_r \in V, \quad \xi:=\xi_{e_1} +\cdots +\xi_{e_r} \in \fa,
\end{equation}
$e$ being a fixed maximal tripotent in $V$.
The root system 
 $\Sig=\Sig(\fh, \fa)$ of 
$(\fh, \fa)$ is of type A, 
 $$
\Sig=
\{ \frac{\beta_j-\beta_k}2
\},
$$
with common multiplicity $a$ or types  B, BC or D,
which we write as
 $$
\Sig=
\{\pm \beta_j,  \frac{\beta_j\pm \beta_k}2,
\pm \frac{\beta_j}2\}
$$
with respective multiplicities $\iota-1, a,  2b$
with the interpretation  that the corresponding multiplicities
$2b=0$ for type C, $\iota-1=0$ for type  B,
and $2b=\iota-1=0$ for type  D.
 Here
$\{\beta_j\}$ is a basis for $\fa^\ast$ normalized
by
$$
\beta_j(\xi_k)=2\delta_{j, k}, \quad j, k=1\cdots, r.
$$
(We write the multiplicity as $\iota-1$ since $\iota=1, 2, 4$ 
has the interpretation as the dimension of the real, complex, quaternionic fields.)
We will view type B as a special case of type BC with the multiplicity
$\iota-1=0$.  We order the roots so that $\beta_1 >\cdots >\beta_r$
(and $\beta_r>0$ for types B and BC)
and denote
\begin{equation}
  \label{eq:half-sum}
\rho=\frac 12 \sum_{\gamma\in \Sig^+}m_{\gamma} \gamma
\end{equation}
the half-sum of positive roots.

We get also a  quadruple characterizing $X$ in $D$ of (2.1),
\begin{equation}
  \label{eq:X-ch}
\boxed{X: (r, \iota-1, a, 2b)  }.
\end{equation}

We list the corresponding quadruples 
$(r, \iota-1, a, 2b)$  and classify them according
to the root system; see  e. g.
\cite{Loos-bsd}, \cite{Ola-Hilg} and \cite{kobayashi-birkh}.

Type $BC\times BC$. The complex domain  is 
$D\times \bar D=(G\times G)/(K\times K)$ 
where $D$
is an irreducible bounded symmetric  
(tube or non-tube type) domain 
of rank $r$ in $\mathbb C^d$,  and 
and $X=H/L=G/K=D$ viewed
as the diagonal part  in 
$D\times \bar D=(G\times G)/(K\times K)$.
(The complex domain  $D\times \bar D$ is reducible so there is
some abuse of definition in \S2.2.)
The
quadruple (\ref{eq:X-ch}) becomes
\begin{equation}
  \label{eq:root-ga}
\boxed{\text{Type} \, \, BC\times BC: 
(r, \iota-1,  a, 2b)=(r^\prime, 1,  a^\prime, 2b^\prime)}
\end{equation}
where $(r^\prime, 1, a^\prime, 2b^\prime)$ is as in
  (\ref{eq:D-ch}) in \S2.1.

Type A. The list of $(G, H)$ is in \cite{FK-book}. 
In this case we have $r^\prime=r$, $a^\prime=a$ and $b=0$, namely
$D$ is of tube type,
\begin{equation}
  \boxed{\text {Type A:}\quad
    (r, a)=(r, 2),  \, (r, 1), \, (r, 4), \, (2, a), \, (3, 8).}
\end{equation}

Type $BC$. $(\fh, \fl)=(\mathfrak{sp}(l, r), \mathfrak{sp}(l)\times \mathfrak{sp}(r))$ ($l>r$)
with $(\fg, \fk)=(\mathfrak{su}(2l, 2r), \mathfrak{su}(2l)\times \mathfrak{su}(2r))$ ($l>r$).
The rank and the root multiplicities  are related by
\begin{equation}
  \label{eq:case3}
  \boxed{\text {Type BC:}\quad
(r, \iota-1, a, 2b)=(r, 3, 4, 4(l-r)),
 \,\, 
(r^\prime, a^\prime, 2b^\prime)=
(2r,  2, 4(l-r))},
\end{equation}
where  $2r=r^\prime$, $2a^\prime =a$.

Type B.  $(H,  L)=(SO_0(l, r),  SO(l)\times SO(r))$ ($l>r$),
$(G, K)=(SU(l, r), S(U(l)\times U(r))$ 
or
  $(H, L)=(SO(2r+1, \mathbb C), SO(2r+1))$,
$(G, K)=(SO^\ast(2(2r+1), U(2r+1))$, and
\begin{equation}
  \label{eq:case2-1}
  \boxed{
\text {Type B-1}: \quad  
(r, \iota-1, a, 2b)=(r^\prime, 1, 1, 1-r), \quad
(r^\prime,  a^\prime, 2b^\prime)=
(r^\prime,  2, 2(l-r))
}.
\end{equation}
or
\begin{equation}
  \label{eq:case2-2}
  \boxed{\text {Type B-2:  } \,\, 
(r, \iota-1, a, 2b)=(r^\prime, 1, 2, 2), \quad
(r^\prime, a^\prime, 2b^\prime)=
(r^\prime, 4, 4)}.
\end{equation}
Here $r^\prime=r$, $a^\prime =2a$.

Type $D$. We have
 $(\fh, \fl)=(\mathfrak{so}(r, r), \mathfrak{so}(r)\times \mathfrak{so}(r))$
with  $(\fg, \fk)=(\mathfrak{su}(r, r), \mathfrak{su}(r)\times \mathfrak{u}(r))$,
or $(\fh, \fl)=(\mathfrak{su}^*(8), \mathfrak{sp}(4))
=(\mathfrak{sl}(4, \mathbb H), \mathfrak{sp}(4))$
with  $(\fg, \fk)=(\mathfrak{e}_{7(-25)}, \mathfrak{e}_{6}\oplus \mathfrak{so}(2))$. 
\begin{equation} 
 \label{eq:caseD-1}
  \boxed{\text {Type D-1:}\,\,  
(r, \iota-1, a, 2b)=(r, 0, 1, 0),
 \,\, 
(r^\prime, a^\prime, 2b^\prime)=
(r, 2, 0)}.
\end{equation}
or
\begin{equation} 
  \label{eq:caseD-2}
  \boxed{\text {Type D-2:}\,\,
(r, \iota-1, a, 2b)=(3, 0, 4, 0),
 \,\, 
(r^\prime,  a^\prime, 2b^\prime)=
(3,  8, 0)}.
\end{equation}
We have here $r^\prime =r$, $a^\prime =2a$.

\begin{rema+}
 The above list can be deduced from Loos \cite{Loos-bsd},
where it is done according to the classification of Jordan triples.
Note that the rank two domain  in the real octonions $\mathbb O^2$
is not listed here,  since it is isomorphic as symmetric
space to the tube domain
of $2\times 2$-quaternionic matrices. The realization of $Sp(2, 2)/Sp(2)\times Sp(2)$
as a real form in $E_{6(-14)}/ Spin(10)\times SO(2)$
is not listed here either as it is realized inside
$SU(4, 4)/ S(U(4)\times U(4))$.
The realization of the exceptional
rank one
 domain $H/L$ with  $(\fh, \fl)=(\mathfrak{f_{4(-20)}}, \mathfrak{spin}(9))$ 
as real form in $G/K$
with 
$(\fg, \fk)=(\mathfrak{e}_{6(-20)}, \mathfrak{spin}(10)\times \mathfrak{so}(2))$
is also not listed as we have excluded the rank one case.
\end{rema+}

\section{$L$-invariant polynomials and their Fock-Fischer norms. 
}

\subsection{Jack polynomials}
Let $\Ome_{\m}=\Ome_{\m}^{(2/a)}$ be the Jack symmetric 
polynomials with multiplicity $\frac a2$ normalized
by
$$
\Ome_{\m}(1^r)=1.
$$
Here we use the abbreviation 
$1^r=(1, \cdots, 1)$.
  In the standard notation \cite{Macd-book}
 it is
$$
\Ome_{\m}(x_1, \cdots, x_r)= \frac{J_m^{(2/a)}(x_1, \cdots, x_r) 
}{J_m^{(\frac 2a)} (1^r)}.
$$
Following \cite{Yan-cjm} we introduce 
\begin{equation}
  \label{eq:q}
q:=q_{a/2}:=1+\frac a2 (r-1).
\end{equation}
and
\begin{equation}
  \label{eq:pi-m}
\pi_{\m}:=\pi_{\m, \frac a2}:=\prod_{1\le i<j\le r}
\frac{m_i-m_j +\frac a2 (j-1)}
{\frac a2 (j-1)}
\frac{(\frac a2 (j-i +1) )_{m_i-m_j}}{(\frac a2(j-i -1)  +1)_{m_i-m_j}}.
\end{equation}
(This is denoted by $d_{\m}$ in  \cite[\S4]{Yan-cjm}.)

\subsection{Type $BC\times BC$}

Consider 
  the complex bounded
symmetric domain $D=H/L=G/K\subset V_{\mathbb C}=\mathbb C^d$
realized as a real form in $D\times \bar D\subset
V_{\mathbb C}\times \bar V_{\mathbb C}$. The parameter
is now $(r, a, b)=(r^\prime, a^\prime, b^\prime)$.
The
space $\mathcal P$ is 
\begin{equation}
  \mathcal P =\mathcal P(V_{\mathbb C})\otimes \overline{
\mathcal P(V_{\mathbb C})}.
\end{equation}
Under $K\times K$ it is decomposed 
as
\begin{equation}
  \label{eq:case-1-poldec}
\mathcal P =\sum_{\n=\m \times\m^\prime}
\mathcal P_{\m}(V_{\mathbb C})\otimes \overline{
\mathcal P_{\m^\prime}(V_{\mathbb C})}.
\end{equation}
The following lemma  follows   immediately 
from Theorem 2.1. All the Pochammer product $(\sig)_{\m}$
in here are understood as $(\sig)_{\m, a/2}$.

\begin{lemm+}
 In the decomposition (\ref{eq:case-1-poldec}),
 $\mathcal P_{\n}^L\ne 0$ if and only if
$\n=(\m, \m)$, 
in which case
the polynomial
$$
p_{(\m, \m)}(x)=\frac {(d/r)_{\m}}{d(\m)} K_{\m}(x, x), \quad
d(\m)=\text{dim}P_{\m},  \quad \n=\m\times \m, 
$$
is the unique $K$-invariant polynomial
in $\mathcal P_{\n}$ normalized by $p_{(\m, \m)}(e)=1$.
The restriction 
of $p_{(\m, \m)}(x)$ on $\fa=\{x=\sum x_j e_j; x_j\in \mathbb R\}\subset V_{\mathbb C}$ is
the Jack symmetric polynomial, $p_{(\m, \m)}(x)
=\Ome_{\m}(x_1^2, \cdots, x_r^2)$.
Its norms in the Fock and Bergman spaces are given by
$$
\Vert p_{(\m, \m)}\Vert_{\mathcal F\otimes \overline{\mathcal F}}^2
=\frac {((d/r)_{\m})^2}{d(\m)}
$$
and
$$
\Vert p_{(\m, \m)}\Vert_{\mathcal H_{\nu}\otimes
\overline{\mathcal H_{\nu}}}^2
=\frac 1{(\nu)_{\m}^2 }
\Vert p_{(\m, \m)}\Vert_{\mathcal F\otimes \overline{\mathcal F}}^2
=\frac {(d/r)_{\m}^2}{(\nu)_{\m}^2 d(\m)}.
$$
 \end{lemm+}

To compare with
Proposition 
\ref{fock-norm} below we write $\Vert p_{(\m, \m)}\Vert_{\mathcal F\otimes \overline{\mathcal F}}^2$
in terms of $\pi_{\m, a^\prime/2}$ defined in (\ref{eq:pi-m}). The dimension
$d(\m)$ is computed in \cite[Lemmas 2.5 and 2.6]{Up-tams} and
is given by
\begin{equation}
\label{up-dim-f}
d(\m)=\frac{
(d/r)_{\m}
}
{
(q)_{\m}
} \pi_{\m}.
\end{equation}
Thus
$$
\Vert p_{(\m, \m)}\Vert_{\mathcal F\otimes \overline{\mathcal F}}^2
=\frac {(d/r)_{\m} (q)_{\m}}{\pi_{\m}}.
$$

\subsection{Type A}
The following is proved in \cite{FK}.
\begin{lemm+}
 In the decomposition (\ref{hua-sch}),
each  component $\mathcal P_{\m}$ has
a unique $L$-invariant polynomial 
$$
p_{\m}(x)=\frac {(d/r)_{\m}}{d(\m)} K_{\m}(x, e), \quad
d(\m)=\text{dim}P_{\m},  \quad \n=\m, 
$$
normalized by $p_{\m}(e)=1$.
The restriction 
of $p_{\m}(x)$ on $\fa$
is the Jack symmetric polynomial, $p_{\m}(x)
=\Ome_{\m}(x_1, \cdots, x_r)$.
Its norms in the Fock and Bergman spaces are given by
$$
\Vert p_{\m}\Vert_{\mathcal F}^2
=\frac {(q)_{\m}}{\pi_{\m}}, \quad
\Vert p_{\m}\Vert_{\mathcal H_{\nu}}^2
=\frac {(q)_{\m}}{ (\nu)_{\m}\pi_{\m}}.
$$
 \end{lemm+}

\subsection{Types B, BC, C, D}

In this section we will generalize the result in \cite{gz-br2}
to the non-tube case; some of which are quite
similar to that of tube domain while others can
be proved by using the results there. We will be rather brief.

The following lemma can be proved by using
the classification theory of spherical pairs
\cite{Kraemer}; for tube domains (namely types A, B, D)
it is also a  consequence of the
Cartan - Helgason theorem
\cite[Chapter V, Theorem 4.1]{He2}.

\begin{lemm+}
\label{cart-helg}  In the decomposition (\ref{hua-sch}),
 $\mathcal P_{\n}^L\ne 0$ if and only if,
\begin{equation}
\label{n-m-rel-h}
\text{Type BC}
\quad \n=(\m, \m):=(m_1, m_1, m_2, m_2,\dots, m_r, m_r)=\sum_{j=1}^rm_j(\ga_{2j-1} +\ga_{2j}),
 \end{equation}
\begin{equation}
\label{n-m-rel-r}
\text{Type B}: 
\n=2\m
=(2m_1, 2m_2, \dots, 2 m_r)=
\sum_{j=1}^r 2m_j\ga_{j},  
\end{equation}
\begin{equation}
\label{n-m-rel-d}
\begin{split}
\text{Type D}: 
\n=2\m +m&=(2m_1+m, 2m_2+m, \dots, 2 m_r +m)\\
&=\sum_{j=1}^r (2m_j+m)\ga_{j},   \quad m=0, 1,
\end{split}
\end{equation}
in which case
$\mathcal P_{\n}^L$ is  one dimensional.
 Here  $m_1\ge m_2\ge \dots \ge m_r\ge 0$.
\end{lemm+}

We will  find the $L$-invariant polynomials in $\mathcal P_{\n}$
in the previous Lemma in terms of the 
Weyl group invariant orthogonal polynomials studied by \cite{Dunkl-cmp-98}, and compute their Fock space norm.
We will denote 
the subspace $\mathbb R e_1 +
\cdots +\mathbb R e_r$ of $V$ also by $\mathfrak a$.

Associated to the root system $\Sig(\fh, \fa)$ there
are the Dunkl difference-differential
operators
\cite{Dunkl-tams},
$$
D_j=\partial_j +\frac 12 \sum_{\a\in \Sig^+}m_{\a}\frac{\a(\xi_j)}{\a(x)}
(1-r_{\a})
$$
acting on polynomials $f(x)$ on $\fa$. 
(It is a realization on the space of polynomials  of the Hecke-algebra of
the tensor product of the symmetric algebra on $\fa$ and the
Weyl group algebra, with $\xi_j$ acting as $D_j$, and $w\in W$
acting by change of variables; see \cite{Opdam-acta}.)

We recall \cite[Proposition 6.2]{gz-br2} an isometric version of the 
Chevalley restriction theorem
(see e.~ g.~ \cite{He2, Torossian-res}
and references therein). 
We define
a norm  \cite{Dunkl-cmp-98}
on 
$\mathcal P(\mathfrak a)^W$
by  
$$
\Vert p\Vert_{B}^2
=p(D) p^\ast\large{|}
_{x=0}
$$ 
for Type $B$  and
$$
\Vert p\Vert_{B}^2=p(\frac 12D) p^\ast\large{|}_{x=0}
$$
for Type $BC$, where for any polynomial $p(x)$, 
$x=x_1 e_1 +\cdots + x_r e_r\in \fa$,
$p(D)$ is obtained by  replacing the linear polynomial $e_j^\ast$
by $D_j$ and $p^\ast$ obtained by
taking the complex conjugate of the coefficients of the
monomials in $e_j^\ast$. (The norm in
\cite{Dunkl-cmp-98} is defined by $p(D) p^\ast\large{|}
_{x=0}$ for root systems of type B or type BC.
The discrepancy here for type BC is
due to the fact that the vectors $e_j$ or the minimal tripotents
in $V$  have norm squares
being twice of that of minimal tripotents in $V_{\bc}$.)

Let
$$
\text{Res}
=\text{Res}_{\mathfrak a}:
\mathcal P(V_{\mathbb C})^L
 \to  \mathcal P(\mathfrak a)^W 
$$
be the restriction map. 

\begin{lemm+}\label{chevalley}
The map $\text{Res}$ is an isometric
 isomorphism between $\mathcal P(V_{\mathbb C})^L$
and the space
$ \mathcal P(\mathfrak a)^W $ of Weyl group invariant polynomials 
on $\fa$. 
\end{lemm+}

It is proved in \cite{Dunkl-cmp-98}
that the polynomials $\Ome_{\m}(x_1^2, \cdots, x_r^2)$ are
then eigenfunctions of the
operators $p(D_1, \cdots, D_r)$,
where $p$ are Weyl group  invariant polynomials on $\fa$. 

\begin{prop+}  \label{res-pn} For each $\m=(m_1, \dots, m_r)$
there exists a unique polynomial $p_{\n}$
in the space $\mathcal P_{\n}^L$  with $\n$ given by $\m$ as in 
Lemma \ref{cart-helg} such that
\begin{equation}
  \label{def-pm-c}
  \text{Res}\, p_{\n}(x_1e_1 + \dots +  x_r e_r)=\Ome_{\m}(x_1^2, \dots, x_r^2).
\end{equation}
for types B and BC, and
\begin{equation}
  \text{Res}\, p_{\n}(x_1e_1 + \dots +  x_r e_r)=(x_1\cdots x_r)^m\, \Ome_{\m}(x_1^2, \dots, x_r^2), \quad m=0, 1,
\end{equation}
for type D.
\end{prop+}

\begin{proof} By Lemma \ref{chevalley}
we have for each $\m$ there exists
a unique $p$  in $\mathcal P(V_{\bc})^L$
such that $\Res  p=\Ome_{\m}(x_1^2, \dots, x_r^2)$. 
We only need to prove that $p$ is the the space
$\mathcal P_{\n}$.
There is linear subspace $V^{0}_{\bc}$ of
 $V_{\bc}$ and
symmetric tube domain   $D_0=G^0/K^0\subset
V^{0}_{\bc}$   in $D=G/K$ of same rank $r^\prime$,
and correspondingly there is a real tube domain
$X_0=H_0/L_0$ of rank $r$ in the domain $X=H/L$; this can
be proved abstractly or by checking the list in our classification.
We view  $\fa$ also as a Cartan subspace for the
symmetric space $X_0$.
The root system of $X_0$ is then of type $D$ or type $C$.
By \cite[Propositions 7.6, 8.3]{gz-br2} we see that
there is a unique $L_0$ invariant polynomials $q$
in $\mathcal P_{\n}(V^{0}_{\bc})$ such that
$\Res{q}=\Ome_{\m}(x_1^2, \dots, x_r^2)
=\Res  p$. That $p$ belongs  $\mathcal P_{\n}(V_{\bc})^L$
follows immediately from the fact the
the isomorphism  $\Res$ is compatible with
the realization of $\a$ in $V_0$ or $V$.
\end{proof}

Next we compute the Fock-Fischer norm using
Proposition \ref{res-pn} and the result of
Dunkl \cite{Dunkl-cmp-98}.

\begin{prop+} \label{fock-norm}With the notation as in Proposition \ref{res-pn} 
we have the following formulas
for the norm squares of the polynomial
$p_{\n}$ in the Fock space and Bergman spaces.

Type $B$:
$$
\Vert p_{\n}\Vert^2_{\mathcal F}
=\frac 1{\pi_{\m} }
2^{2|\m|} 
{
(q)_{\m, \frac a2} 
((r-1)\frac a2 +b+\frac 12)_{\m, \frac a2}
}
$$
and 
$$
\Vert p_{\n}\Vert^2_{\nu}
=
\frac 1{\pi_{\m} }
\frac {
(q)_{\m, \frac a2} 
((r-1)\frac a2 +b+\frac 12)_{\m, \frac a2}
}
{
(\frac{\nu}2)_{\m, \frac a2}
 (\frac{\nu+1}2)_{\m, \frac a2}
}; 
$$

Type $BC:$
$$
\Vert p_{\n}\Vert^2_{\mathcal F}
=\frac 1{\pi_{\m} }  {(q)_{\m, \frac a2} 
((r-1)\frac a2 +\frac {\iota +2b}2)_{\m, \frac a2}},
$$
and
$$
\Vert p_{\n}\Vert^2_{\nu}
=\frac 1{\pi_{\m} }
\frac
{
(q)_{\m, \frac a2} 
((r-1)\frac a2 +\frac {\iota +2b}2)_{\m, \frac a2}
}
{
(\nu)_{\m, \frac a2}
 (\nu-\frac{ a^\prime} 2)_{\m, \frac a2}
 };
$$
Type $D:$
$$
\Vert p_{2\m}\Vert^2_{\mathcal F}
=\frac 1{\pi_{\m} }
2^{2|\m|} 
{
(q)_{\m, \frac a2} 
((r-1)\frac a2 +b+\frac 12)_{\m, \frac a2}
},
$$
$$
 \Vert p_{2\m+1}\Vert^2_{\mathcal F}
=\frac
{ 1}
{\pi_{\m} }
2^r \prod_{j=1}^r(\frac a2 (r-1) +\frac 12- \frac a2 (j-1) )
2^{2|\m|} 
(q)_{\m, \frac a2}  (\frac a2 (r-1) +\frac 32)_{\m},
$$
$$
\Vert p_{2\m}\Vert^2_{\nu}
=
\frac 1{\pi_{\m} }
\frac {
(q)_{\m, \frac a2} 
((r-1)\frac a2 +\frac 12)_{\m, \frac a2}
}
{
(\frac{\nu}2)_{\m, \frac a2}
 (\frac{\nu+1}2)_{\m, \frac a2}
}; 
$$
$$
\Vert p_{2\m+1}\Vert^2_{\nu}
=
\frac
{ 1}
{\pi_{\m} }
\prod_{j=1}^r
\frac
 {\frac a2 (r-1)+\frac 12 -\frac a2(j-1)}
{\frac \nu 2 -\frac a2(j-1)}
\frac 
{
(q)_{\m, \frac a2} 
((r-1)\frac a2 +\frac 32)_{\m, \frac a2}
}
{
(\frac{\nu}2 +\frac 12)_{\m, \frac a2}
 (\frac{\nu}2 +1)_{\m, \frac a2}
}; 
$$
\end{prop+}

\begin{proof}
By Lemma \ref{chevalley}
and Proposition \ref{res-pn}
 we have
$$
\Vert p_{\n}\Vert^2_{\mathcal F}=
\Vert \text{Res}\, p_{\n}\Vert^2_{B}
=\Vert \Ome_{\m}\Vert^2_{B}
$$ 
for type $B$,
and
$$
\Vert p_{\n}\Vert^2_{\mathcal F}=
2^{-|\n|}\Vert \Ome_{\m}\Vert^2_{B}
=2^{-2|\m|}\Vert \Ome_{\m}\Vert^2_{B}
$$ 
for type  $BC$, since 
$p_{\n}$ is a polynomial of degree $|\n|=2|\m|$. 
The right hand side is computed by
Dunkl \cite[Section 5]{Dunkl-cmp-98}, where
the norm for type BC is defined as
for type B.
 However
to express the resulting formula as stated requires
some rather technical computations, so we will
adapt the notations there,   his $N, k, k_1,  \lambda$
are 
our $r, \frac a2, \frac{2b+\iota-1}2, \m$ etc., and the $h$ below is
the shifted hook length product defined in
Definition 3.17 there.  (Observe that
for root system of type BC the constant
$k_1$ in  \cite[Section 5, (5.1)]{Dunkl-cmp-98}
is a sum of two multiplicities.)
The Jack polynomial
$J_{\m}$ is expressed in terms
of the polynomial $j_{\m}$ by \cite[Section 3, p.465]{Dunkl-cmp-98}
$$
J_{\m}
=\frac
{
(r\frac a2 +1)_{\m} h(\m, 1)
 }
{(r\frac a2 +1)_{\m} (\frac a2)^{|\m|}(\# S_r \m)
}
j_{\m}
$$
with
$$
J_{\m}(1^r)=(r \frac a2 )_{\m}(\frac a2)^{-|\m|}.
$$
The norm $j_{\m}$ is (see p. 480-495, loc. cit.)
$$
\Vert j_{\m}\Vert_B^2 = 2^{2|\m|} (r\frac a2 +1)_{\m}
((r-1)\frac a2 +\frac{\iota-1 +2b}2 +\frac 12)_{\m} 
(\# S_r \m)
\mathcal E (\m^R) 
\frac{h(\m, \frac a2 +1)}
{h(\m, 1)}.
$$
We find then
$$
\Vert\Omega\Vert_B^2=
2^{2|\m|} 
\frac
{
((r-1)\frac a2 +\frac{\iota-1 +2b}2 +\frac 12)_{\m}
}
{
(r\frac a2)_{\m}
}
h(\m, 1) h(\m, \frac a2)
$$
The shifted 
 hook length product $h$
is related to the 
 upper and lower hook length products $ h^*(\m)$
and $h_*(\m)$  
(Stanley \cite{Stanley-jack}, Macdonald \cite{Macd-book})
by
$$
h(\m, 1)h(\m, \frac a2)=(\frac{a}{2})^{2|\m|}
 (h_*(\m) h^*(\m))
$$
and 
which, by
 \cite[Proposition 4.1]{Yan-cjm},  can be further written in terms of the 
quantity $\pi_{\m}$ defined in (\ref{eq:pi-m}),
$$
 h_*(\m) h^*(\m)=(\frac{2}{a})^{2|\m|} 
\frac
{
(q)_{\m} (r\frac a2)_{\m}
}
{
(\pi)_{\m}}.
$$  
Namely
$$
\Vert\Omega\Vert_B^2=2^{2|\m|}
\frac {
(q)_{\m}
}
{
\pi_{\m}
}
((r-1)\frac a2 +\frac{\iota-1 +2b}2+\frac 12)_{\m},
$$
proving our claim for the Fock space norm.

To find the Bergman space norm 
of $p_{\n}$ we use Theorem 2.1.  For Type $B$  
we have,  $a^\prime =2a$, $\iota=1$,   $\n=2\m$,
\begin{equation*}
\begin{split}
\Vert p_{\n}\Vert^2_{\nu}
&=\frac 1{(\nu)_{\n,  a^\prime/2}}
\Vert p_{\n}\Vert^2_{\mathcal F}\\
&= \frac 1 {(\frac{\nu}2)_{\m, \frac a2}
 (\frac{\nu+1}2)_{\m, \frac a2} }
\frac 1{
(q)_{\m}
}
{
\pi_{\m}
}
((r-1)\frac a2 +b+\frac 12)_{\m}.
\end{split}
\end{equation*}
since
\begin{equation}
\label{poch-r}
(\nu)_{\n,  a^\prime/2} =2^{2|\m|}
(\frac \nu 2)_{\m, a/2}
(\frac {\nu+1} 2 )_{\m, a/2},
\end{equation}
for any $\nu$.

For Type $BC$, it holds  $a =2a^\prime$, and
\begin{equation}
\label{poch-h}
(\nu)_{\n,  a^\prime/2} =
(\nu)_{\m,  a/2} (\nu-a^\prime /2)_{\m,  a/2},
\end{equation}
and we get  the last equality.
\end{proof}

\section{Principal series representations
on maximal boundaries and intertwining
 operators into the space $\mathcal H_\nu$. Type $BC\times BC$
}

The domain $X=H/L=G/K$ is the complex bounded symmetric
of rank $r$ with root multiplicity
$(1, a, 2b)$. 
We shall study spherical representations
defined on the Shilov boundary of $D$ in terms
of Macdonald ${}_2F_1$ hypergeometric functions
and their realization in the holomorphic
representation $\HtH$. First we recall some
general factors about the hypergeometric
series which we will need through the rest of the paper.

\subsection{Hypergeometric series
with general parameter $\frac 2a$}

In the following sections
 we will need the hypergeometric
functions 
defined in terms of
the Jack symmetric polynomials.
Let $\alpha =(\alpha_l, \cdots, \alpha_{k})$,
 $\beta =(\beta_l, \cdots, \beta_{l})$ be two tuples
of real positive numbers, such that
that $\beta_1, \cdots, \beta_l >\frac a2(j-1)$
and  $t=(t_1, \cdots, t_r) \in [0, 1)^r$. We define
$$
{}_{k}F_l^{(2/a)}(\alpha; \beta;  t)
=\sum_{\m}\frac{(\alpha_1)_{\m} \cdots
(\alpha_{k})_{\m}}
{(\beta_1)_{\m}\cdots (\beta_l)_{\m}
} \frac{\pi_{\m} }{(q)_{\m}}
\Ome_{\m}(t_1^2,
\cdots, t_r^2).
$$
Note that in case one of $\alpha$'s  is $\frac a 2(j-1)$
for some $1\le j\le r$
we have $(\alpha_1)_{\m} \cdots
(\alpha_{k})_{\m}=0$ and the sum is only
over those with $m_j=0$,  namely $\m=(m_1, \cdots, m_{j-1}, 0, \cdots, 0)$. 
We will suppress the upper index $2/a$ when no confusion would arise.
The series ${}_{2}F_1$ 
has been well-studied as it is related
the spherical functions on 
symmetric domains (see below).

The convergence property
is similar to that of ${}_{2}F_1$, namely we
have
\begin{lemm+}\label{conv-le}
Suppose $\a_1, \cdots, \a_{l+1} >0$
and $\beta_1, \cdots, \beta_l >\frac a2(j-1)$.
The hypergeometric series 
${}_{l+1}F_l(\alpha; \beta;  t)
$ is bounded on the set $[0, 1)^r$ if  and only if
$$ \sum_{p=1}^{l+1} \alpha_p-
\sum_{p=1}^{l}\beta_p  <-\frac a2 (r-1)
$$
in which case the series $F(\alpha; \beta; 1^r)$
 is convergent.
\end{lemm+} 
\begin{proof}
The sufficiency for $ l=1$ was  proved in \cite{FK}
for special values of $(r, l)$ (corresponding
to a complex bounded symmetric domain)
and was generalized
by Yan  \cite{Yan-cjm} 
to general $(k, r)$, in the case $\alpha=(\alpha_1, \alpha_2)$
and $\beta=\beta_1$; the general case 
is exact the same (see e.g. \cite{EZ-fourier}). The necessary part
is essentially also proved 
in  \cite{Yan-cjm}  with approximate behavior
of the function near the certain boundary  part
of  $[0, 1)^r$, 
and we provide
here the argument. Put
$\varepsi:=\sum_{p=1}^{l+1} \alpha_p-
\sum_{p=1}^{l}\beta_p +\frac a2(r-1)$.
 Consider $t$ along  the diagonal  
$(t_1, \cdots, t_r)=(t, \cdots, t)$, we have, by 
the Stirling formula (see \cite{EZ-fourier}, (2.9))
$$
{}_{l+1}F_l(\alpha; \beta;  t)
\asymp \sum_{m}\prod_{j=1}^r (1+m_j)^{\varepsi -q}
\prod_{1\le i <j\le r}(1+m_i-m_j)^{a} t^{|\m|}.
$$
Using the evaluation formula for $j_m$ (see e.g. formula (2.13) in \cite{EZ-fourier})
we see by elementary computations \cite{Yan-cjm}
that
the above series behaves the same as  $\log\frac{1}{1-t}$ for 
$\varepsi=0$ 
or $(1-t)^{\varepsi + \frac a2(r-1)}$
for $\varepsi
>0$, which is thus unbounded.
\end{proof}

For the ${}_2F_1$-series
the sum 
 $F(\alpha; \beta; 1^r)$ has also been explicitly evaluated; see
\cite{Opdam-invmath89},
 \cite{Yan-cjm} and \cite{Beerends-Opdam}.

\subsection{Induced representation of $H$ on $L^2(S)$
and spherical functions. }

Fix the maximal tripotent $e\in V_{\mathbb c}$
and  let $S=K\cdot e$. 
It is well understood
that $S$ is the Shilov boundary of $D$ and $S=K/K_e=G/P$ where $K_e$ and 
$P$ is the isotropic subgroup of $K$
and respectively $G$ of $e$.
Consider
the root space decomposition of $\fg$ under 
the element  $\xi \in \fa$ in  (\ref{eq:max-tri}),
\begin{equation}
  \label{eq:parab-P-1}
\fg=\fn_{-} + \fn_0 + \fn_+.  
\end{equation}
Then $P=MAN$ is a parabolic subgroup with
 Lie algebra $ \fn_0 + \fn_+$
and $A$ is the Lie group with Lie algebra
$\mathbb R\xi$ and $MA$ is the Levi component
with Lie algebra  $\fn_0$.

For $\lam\in \mathbb C$ identified with the
linear function $\lam \xi^\ast$
 we let
 $Ind_P^{G}(\lam)=L^2(S)$
be the induced representation of $G$ 
on $L^2(S)=L^2(S, dv)$, where $dv$
is the  normalized $K$-invariant measure on $S$;
see \cite[Chapter VII, \S\S1-2]{Kn-book}.
The group action is given
by
\begin{equation}
U(\lam, g)f(v)=|J_{g^{-1}}(v)|^{\frac {2n-i\lam}{pr}} f(g^{-1}v)
\end{equation}
where $J$ is as in (\ref{pi-nu})
the complex Jacobian at $v\in S$. In particular $Ind_P^G(\lam)$
is unitary when $\lam$ is real.

It is known that  Harish-Chandra
$e$-function with respect to the decomposition
$G=KMAN$ is
given by
$$
e^{(i\lam +\rho)(A(kg))}
=\frac{h(z, \bar z)^{\frac{\sig}2}}{h(z, \bar v)^\sig}, \quad v=k^{-1}e\in S, \, k\in K
$$
where  $\sig$  is determined by $\lam$ 
and vice versa  via
\begin{equation}
  \label{eq:lam-sig-c}
\sig=\frac 1{2r}(i\lam +\rho)(\xi)=\frac { i\lam}{2r} +\frac 12(1+b+ \frac a2(r-1)),
\quad \, i\lam=2r\sig-\rho(\xi),
\end{equation}
and $e^{A(kg)}$ stands for the $A=\exp(\mathbb R\xi)$-component in the decomposition; see e.g. \cite{kz-hua}. 
(Similar formulas
holds,  \cite{UU}, in the Siegel domain realization,
 for general linear functional $\lam $ on $\fa$.) 
The Poisson transform from $
Ind_{P}^G(\lam)
$ into the space of eigenfunctions
of $G$-invariant differential operators 
is given by
  \begin{equation}
    \label{eq:poisson-1}
P_{\lam}f(z)=\int_S\left (\frac{h(z, \bar z)}{|h(z, \bar v)|^{2}}\right)^{\sigma} f(v) dv,
  \end{equation}
which intertwines the 
induced representation  $Ind_P^{G}(\lam)$
and the regular action on $X$.

The corresponding spherical function can be expressed in
terms of the hypergeometric function. We recall
the following know result; see
e.~g.~ \cite{FK, Yan-cjm, EZ-fourier}. 

\begin{lemm+} The spherical function
$$
\phi_{\lam}(z)
=P_{\lam}1(z)=
\int_S\frac{h(z, \bar z)^{\sig}}
{h(z, \bar v)^\sig h(v, \bar z)^\sig} d\nu(v),
$$
when restricted to the radial directions
$z=t_1 e_1 +\cdots +t_r e_r$, is given by
$$
\phi_{\lam}(z)=(\prod_{j=1}^r (1-t_j^2))^{\sig}
\,{}_2F_1(\sig, \sig; 1+b +\frac a2(r-1);  t), \quad
\sig =i\frac {\lam}{2r}
+ \frac 12 (1+b +\frac a 2(r-1)) 
$$
\end{lemm+}

\subsection{Discrete components of $(\mathcal
H_\nu \otimes\overline{\mathcal H_{\nu}}, G\times G)$
under $G$ 
for $\nu>\frac a2 (r-1)$ for type one domains $D$}

In this section we will prove that
the spherical representation
in the previous section can be realized
as a discrete component in the tensor
product decomposition of $\mathcal H_\nu \otimes\overline{\mathcal H_{\nu}}$
for $\nu$ in the continuous part of the Wallach set
only if $D$ is Type I domains,  and in that
case we will find the exact such parameters $\lam$
and the explicit intertwining operator.

We consider  the sesqui-holomorphic extension
of the Poisson transform. More precisely,
 for $\nu, k\in \mathbb C$ 
and $f\in L^2(S)$ we define
$T_{\nu, k}f(z, \bar w)$ to be the holomorphic
function
\begin{equation}
  \label{eq:T-nu-k}
T_{\nu, k}f(z, \bar w)= h(z, \bar w)^{k}
\int_S\left (\frac{1}{h(z, \bar v) h(v, \bar w)}\right)^{\nu+k} f(v) dv, \quad (z, w)\in  d\times \bar D.
\end{equation}
Its restriction to the diagonal is up to a factor the Poission transform, viz
$$
T_{\nu, k}f (z, \bar z)= h(z, \bar z)^{-\nu} P_{\lam} f(z)
$$
with $i\lam$ determined by $\sig=\nu+k$ as in (\ref{eq:lam-sig-c}):
\begin{equation}
\label{lam-nu-k-1}
 i\lam=2r(\nu+k)-\rho(\xi).
\end{equation}

The intertwining property of $P_{\lam}$ and the transformation
formula of $h(z, \bar w)$ under $G$ imply immediately

\begin{lemm+} Let $\nu, k\in \mathbb C$
and $\lam$ be as in 
(\ref{lam-nu-k-1}).
 The operator
$
T_{\nu, k}:  Ind_P^G(\lam)=L^2(S) \to \mathcal O(D\times \bar D),
$
is a formal $G$-intertwining operator
from the induced representation
$Ind_P^G(\lam)$ to the
space $\mathcal O(D\times \bar D)$
with the action $\pi_\nu \otimes \overline {\pi_\nu}$.
\end{lemm+}

We determine  when the image of $T_{\nu, k}$ is
in the Hilbert space
 $\mathcal H_{\nu}\otimes \overline {\mathcal
H_{\nu}}$.

\begin{lemm+}Let $\nu>\frac a2(r-1)$, $k\ge 0$ be an integer
and  $\lam$ be given in (\ref{lam-nu-k-1}).
The image of the constant function $1$ in $Ind_P^G(\lam)$
is mapped under $T_{\nu, k}
$ into the Hilbert space $\mathcal H_{\nu}\otimes \overline {\mathcal
H_{\nu}}$
if and only if 
$$2\nu +4k <1+b.
$$
Furthermore the inequality has  a possible solution
if and only if  when $D$ is a type I non-tube type domain $SU(l, r)/S(U(l)\times U(r))$
with $a=2$,
$b=l-r>2 $
 and 
$$
r-1<\nu +2k<\frac 12(1+l-r)
.
$$
\end{lemm+}
\begin{proof} We write 
$$
T_{\nu, k}f(z, \bar w)= h(z, \bar w)^{k}
F(z, \bar w), \quad F(z, \bar w):=
\int_S\left (\frac{1}{h(z, \bar v) h(v, \bar w)}\right)^{\nu+k} f(v) dv, 
$$
and we  shall prove that $F(z, \bar w)$ is in the space, and our results
then follows since the function $h(z, \bar w)^k$ for non negative
integer $k$
is a polynomial in $z$ and $\bar w$, the multiplication operator
by coordinate functions  is a bounded operator on $\mathcal H_{\nu}$ (see \cite{az-mult}).
The function $F$ can be computed  by using the expansion
(2.7) (see also \cite{FK})
$$
F(z, \bar w)=\sum_{\m}
\frac{(\nu+k)_{\m}^2}{(d/r)_{\m}} K_{\m}(z, \bar w).
$$
Its norm square in $\HtH$, again by Theorem 2.1,  is
$$
\Vert F\Vert_{\HtH}^2=\sum_{\m}
\left(\frac{(\nu+k)_{\m}^2}{(d/r)_{\m}}\right)^2 
\frac{1}{(\nu)_{\m}^2 } d(\m)
$$
which by (\ref{up-dim-f}) and the definition
of hypergeometric function, is
$$
{}_4F_3(\nu+k, \nu+k, \nu +k, \nu +k; \nu, \nu,  d/r; 1^r)
$$
By Lemma 
\ref{conv-le},  the series is convergent if and only if
$$
4(\nu +k)-2\nu -\frac dr <-\frac a2(r-1).
$$
Simplifying this is $2\nu +4k < 1+b$, proving
the first part.
Checking over the list  (\ref{eq:D-cl}) of bounded symmetric domains we
see that this has a solution only if $D$ is type one non-tube domain.
The condition for $\nu$ and $k$ is then
$r-1 <\nu \le \nu +2k <
\frac 12(1+l-r)$. 
\end{proof}

\begin{theo+} Let $D$ be 
 the type one domain
$SU(l, r)/S(U(l)\times U(r))$
with $l-r>2$  and  
$$r-1<\nu<\frac 12(1+l-r +2(r-1).$$
Let  $k$ be a nonnegative positive integer such that
$$
0\le k <\frac 14(1+l-r  -2\nu).
$$
Then the spherical function $\phi_{\lam}$, for 
$$i\lam=2r(\nu +k)-\rho(\xi),$$
is positive definite and the corresponding unitary
spherical
representation of $G$ appears
as a discrete component in the irreducible
representation of $(\HtH, G\times G)$.
\end{theo+}
\begin{proof}
Consider the linear span of the constant function $1$
in $L^2(S)$ under the principal series representation,
$$
\mathcal S_{\lam}=\span\{U(g, \lam) 1; g\in G\}
$$
It is a pre-Hilbert space with
the inner product
$$
(f, g):= (T_{\nu, k} f,  T_{\nu, k} f)_{\HtH}
$$
It follows from the previous two Lemmas  that
this is well-defined and  $G$-invariant. Its completion
is then a spherical unitary representation, which
in turn is irreducible since it is defined by a spherical
function,  and is realized as
a discrete component in the tensor product
via $T_{\nu, k}$. 
\end{proof}

\subsection{Discrete component of $(\mathcal H_\nu \otimes\overline{H_{\nu}}, 
G\times G) $ under $G$
for $\nu=\frac a2(j-1)$ being  a singular Wallach point}

We consider the tensor product
$\mathcal H_\nu \otimes\overline{\mathcal H_{\nu}}$
with 
$\nu=\frac {a}2
 (j-1)$ being a singular
Wallach point. 

The  operator
$T_{\nu, k}$ intertwines the induced representation
with the action $\pi_{\nu}\otimes \overline{\pi_{\nu}}$
on $\mathcal O(D\times \bar D)$. However
the space $\HtH$ has $K$-types restriction
 (\ref{eq:sig-wa-hil}). Thus only the operator
$T_{\nu}=T_{\nu, 0}$,
\begin{equation}
  \label{eq:T-nu-0}
T_{\nu}f(z, w):=\int_S \frac 1{h(z, \bar v)^{\nu} h(v, \bar w)^{\nu } } dv,  
\end{equation}
 will be possibly
an operator into $\HtH$.
Furthermore by some similar computation
as in Lemma 4.3 (see also the proof below)
this will happen possibly only 
for type I domains.

\begin{theo+} Let $D$ be 
 the type one domain
$SU(l, r)/S(U(l)\times U(r))$, $l-r>2$
and let $\nu= j-1$, $2\le j\le r$,  be a singular Wallach point.
Suppose $l-r>2j-3$. The spherical function $\phi_{\lam}$,
for $\lam$ given by   $i\lam =2r(j-1) -r(1+l-r +(r-1))$, 
is positive definite and  the corresponding unitary
spherical representation
appears
as a discrete component in the irreducible
representation of $\HtH$.
\end{theo+}

\begin{proof} We prove that the image of the function $1$ under
$T_{\nu}$ is in the Hilbert space
$\HtH$, and the rest is proved by similar arguments
as that  of the previous Theorem.
We have, 
$$
F(z, w):=
(T_{\nu} 1)(z, w)=\sum_{\m; m_j=0}
\frac{(\nu+k)_{\m}^2}{(d/r^\prime)_{\m}} K_{\m}(z, w).
$$
Its norm in
$$
\sum_{\m; m_j=0}
\frac{(\nu)_{\m}^2} {(d/r)^2 _{\m}}
d(\m)
=\sum_{\m; m_j=0}
\frac{(\nu)_{\m}^2} {(d/r)_{\m}}
\frac{\pi_{\m}}{(q)_{\m}}
$$
which by Lemma \ref{conv-le} is convergent if
$$
2\nu -d/r=2v-(1+b +(r-1)) <-(r-1)
$$
namely if $l-r=b >2\nu -1 =2(j-1) -1=2j-3.$
\end{proof}

\section{Principal series representations
on maximal boundaries and spherical functions.
Intertwining
 operators into the space $\mathcal H_\nu$.
 Types A, B, BC,  D}

\subsection{Induced representations
and spherical functions.}

For the real domain $X$ in $D$ 
we let $S= L\cdot e$ be the orbit of $e$ under $L$.
Then $L \subset \partial_e D \cap V$, where
 $\partial_e D$ is the Shilov boundary of $D$. (In certain 
cases $S$ is a true subset of $\partial_e D \cap V$.) 
Then $S$ can be realized as $S=G/P$ 
where $P$ is a parabolic subgroup of $H$ 
with Lie algebra given by $ \fn_0 + \fn_+$
in the decomposition 
\begin{equation}
  \label{eq:parab-P-2}
\fh=\fn_{-} + \fn_0 + \fn_+,
\end{equation}
under the adjoint action $\xi=\xi_1+\cdots +\xi_r$.

Let $\lam=\lam\xi^\ast$  on $\bc \xi$. 
We consider the induced representation $Ind_P^H(\lam)$ with
of $H$ realized on $L^2(S)$.

\begin{theo+}
The spherical function 
$\phi_{\lam}(z)$ is given by the integral
$$
\phi_{\lam}(z)
=\int_S\frac{h(z, \bar z)^{\frac{\sig}2}}{h(z, \bar v)^\sig} dv.
$$
Its restriction on  $z=t_1 e_1 +\cdots + t_r e_r$,
$|t_1|, \cdots, |t_j| < 1$, 
 is further given (and uniquely determined) by 
\begin{itemize}
\item[Type A:]
\begin{equation*}
 \phi_{\lam}(z) = \prod_{j=1}^r (1-t_j)^{ \frac \sig 2}
\,{}_1 F_1(\sig; 1+\frac a2(r-1); t),
\, \sig = \frac{i\lam}{r} + \frac  a2(r-1);
\end{equation*}
\item[Type B:]
\begin{equation*}
 \phi_{\lam}(z) = (\prod_{j=1}^r(1-t_j^2)^{\frac \sig 2}
\,{}_2 F_1(\frac \sig 2, \frac{\sig +1}2; \frac a2(r-1)+b +\frac 12; t^2),
\, \sig = \frac{i\lam}{r} + b +\frac  a2(r-1);
\end{equation*}
\item[Type BC:]
\begin{equation*}
 \phi_{\lam}(z) =
 (\prod_{j=1}^r(1-t_j^2)^{\sig}
\,{}_2 F_1 (\sig, \sig -1; \frac a2 (r-1) +\frac{\iota +2b}2; t^2),
\, 
\sig= \frac{i\lam}{2r} +\frac{1}2 (\iota-1 +b +2a(r-1));
\end{equation*}
\item[Type D:]
\begin{equation*}
\begin{split}
 \phi_{\lam}(z) & = 
 (\prod_{j=1}^r(1-t_j^2)^{\frac{\sig}2}
\,{}_2 F_1(\frac \sig 2, \frac{\sig +1}2; \frac a2(r-1) +\frac 12; t^2)\\
&\quad +\prod_{j=1}^r
\frac
{ \frac {\sig}2 -\frac a2(j-1) }
{\frac a2 (r-1) -\frac a2(j-1) +\frac 12
 }
 (\prod_{j=1}^r t_j(1-t_j^2)^{\sig/2})\,
{}_2 F_1(\frac \sig 2 +1, \frac{\sig +1}2; \frac a2(r-1) +\frac 32; t^2),
\end{split}
\end{equation*}
\begin{equation*}
\sig=   \frac{i\lam}{r}  +\frac  a2(r-1).
\end{equation*}
\end{itemize}
\end{theo+}

\begin{proof} We claim first that the Harish-Chandra
$e$-function is
given by
\begin{equation}
  \label{eq:h-c-e}
e^{(i\lam +\rho)(A(kg))}
=\frac{h(z, \bar z)^{\frac{\sig}2}}{h(z, \bar v)^\sig}.  
\end{equation}
This formula, in the Siegel domain realization
of $X$,  is given in \cite{gkz-bere-rbsd} generalizing
that of Upmeier-Unterberger \cite{UU}
for the complex case. Here it can be simply proved
by using the transformation rule of $h(z, \bar w)$
under the group $H$. We get thus the integral representation
of $\phi_{\lam}$.
We compute the integration using the Faraut-Koranyi expansion (\ref{FK-exp}).
We have (as $\bar z=z, \bar v=v$ for  $z\in X$, $v\in S$ we will drop the bar)
$$
\phi_{\lam}(z)
={h(z, z)^{\frac{\sig}2}}
\int_S \sum_{\n}(\sig)_{\n, \frac{a^\prime}2}
K_{\n}(z, v) dv
={h(z, z)^{\frac{\sig}2}}
 \sum_{\n}(\sig)_{\n, \frac{a^\prime}2}
\int_S K_{\n}(z, v) dv,
$$
where the interchanging of the integration and the summation
is justified by the uniform convergence of the expansion
(\ref{FK-exp})
 on $S$ for fixed $z\in D$.
By the $K$-invariance of $K_{\n}$, $K_{\n}(kz, kv)=K_{\n}(z, v)$,
 and the $L$-invariance of the measure $dv$
we have $\int_S K_{\n}(z, v)dv$ is a $L$-invariant
polynomial in $\mathcal P_{\n}$. Thus
$$\int_S K_{\n}(z, v)dv= C_{\n} p_{\n}(z), \quad z\in D
$$
for those $\n$ given in Lemma 3.3 in terms of $\m$
and for some constant  $C_{\n}$; otherwise 
it is zero.
 We claim that
$$
C_{\n}=\frac{1}{\langle p_{\n}, p_{\n}\rangle_{\mathcal F}}.
$$
Indeed, the left hand side is an $L$-invariant
element in $\mathcal P_{\n}$
thus is a multiple of $p_{\n}$ determined in Proposition 3.5. To find the
constant we compute the norm square of the left hand
in the Fock space.
By definition of $K_{\n}$ we have it is
$$
\int_S \int_S 
K_{\n}(w, v)dv\,dw,
$$
which is, by the invariant of $K_n$ under $L\subset K$
and that $S=L\cdot e$,
$$
\int_S K_{\n}(w, e)dw =C_{\n} p_{\n}(e) =C_{\n}.
$$
On the other hand, the squared norm of the right hand
side
is
$$
C_n^2 \langle p_{\n}, p_{\n}\rangle_{\mathcal F},
$$
proving our claim.
Thus
$$
\phi_{\lam}(z)
={h(z, z)^{\frac{\sig}2}}\sum_{\n} (\sig)_{\n, \frac{a^\prime}2}
\frac{1}{\langle p_{\n}, p_{\n}\rangle_{\mathcal F}} 
p_n(z).
$$

For Type A we have $\n=\m$,
  $a=a^\prime$ and
$$
\frac{1}{\langle p_{\n}, p_{\n}\rangle_{\mathcal F}} 
=\frac{\pi_{\m}}{(q)_{\m}}
$$
by Lemma 3.2.
The function $\phi_{\lam}(z)$ is, for $z=t_1 e_1 +\cdots + t_r e_r$,
$$
\phi_{\lam}(z)
= \prod_{j=1}^r (1-t_j^2)^{ \frac \sig 2}
\sum_{\m} (\sig)_{\n, \frac{a^\prime}2}
\frac{\pi_{\m}}{(q)_{\m}}
\Ome_{\m}(t)
= \prod_{j=1}^r (1-t_j^2)^{ \frac \sig 2}
{}_1F_1(\sig; q; t)
$$

For Type B,  we have $\n=2\m$, $a^\prime =2a$, 
the polynomial $p_{\n}$ is the Jack polynomial
$\Ome_{\m}$ by Proposition 3.5, and its norm square
 $
{\langle p_{\n}, p_{\n}\rangle_{\mathcal F}}
$ is computed in  Proposition \ref{fock-norm}.
This gives, for $z=t_1 e_1 +\cdots + t_r e_r$,
\begin{equation*}
  \begin{split}
&\quad \phi_{\lam}(z)\\
&={h(z, z)^{\frac{\sig}2}}
 \sum_{\n=2\m}(\sig)_{\n, \frac{a^\prime}2}
\left( \frac 1{\pi_{\m} }
2^{2|\m|} 
{
(q)_{\m, \frac a2} 
((r-1)\frac a2 +b+\frac 12)_{\m, \frac a2}
}\right)^{-1} \Ome_{\m}(t^2)\\
&= {h(z, z)^{\frac{\sig}2}}
 \sum_{\m} 
2^{2|\m|} (\frac{\sig}2)_{\m, \frac{a}2}  (\frac{\sig+1}2)_{\m, \frac{a}2}
\left( \frac 1{\pi_{\m} }
2^{2|\m|} 
{
(q)_{\m, \frac a2} 
((r-1)\frac a2 +b+\frac 12)_{\m, \frac a2}
}\right)^{-1} \Ome_{\m}(t^2)
\\
&={h(z, z)^{\frac{\sig}2}}
 \sum_{\m} \frac{(\frac{\sig}2)_{\m}  (\frac{\sig+1}2)_{\m} }
{((r-1)\frac a2 +b+\frac 12)_{\m}
}
\frac{\pi_{\m}}
{(q)_{\m}}
\Ome_{\m}(t^2)\\
&={h(z, z)^{\frac{\sig}2}}
{}_2 F_1 (\frac{\sig}2, \frac{\sig+1}2;
(r-1)\frac a2 +b+\frac 12; t^2)
  \end{split}
\end{equation*}
where in the second step we have written
$(\sig)_{\n, \frac{a^\prime}2}$
in terms of  $(c)_{\n, \frac{a}2}$
as in  (\ref{poch-r}).

The remaining types BC or  D are done
by the same method.
Note that for type BC we have $r^\prime =2r$
and $h(z, z)= \prod_{j=1}^r (1-t_j^2)^{2}$.
\end{proof}

\subsection{Discrete components of $(\mathcal H_\nu, G)$ under $H$
for $\nu>\frac {a^\prime}2 (r^\prime-1)$}

In this section we will find and realize certain discrete
components in the branching of the holomorphic
representations $\mathcal H_{\nu}$ 
of $G$ under $H$ 
using the Poisson transform studied above.
Similar computations as in Lemma 4.4  show
that the 
operator $T_{\nu, k}$
defined in  (\ref{eq:T-nu-0}) below
maps the spherical representation
into the holomorphic representation
 $\mathcal H_{\nu}$  only if $X$ is of
Type $B$ with $H=SO(r, l)$,  $l>r$,
or   Type $BC$ with $H=Sp(r, l)$, $l>r$.
 We will thus only consider those cases.
The corresponding group $G$
is then $SU(l, r)$ ($r^\prime=r$) or $SU(2l, 2r)$ ($r^\prime=2r$).

\begin{theo+} Let $H$ be 
the group $SO_0(l, r)$
or $Sp(l, r)$.
Let $l, r$ satisfy $l-r>2(r -1)$ for $H=SO_0(l, r)$
and $l-r\ge 2(r-1)$ for    $H=Sp_0(l, r)$.
Suppose $\nu>r^\prime-1$ is a 
be a point in the continuous part of the Wallach set
of $G$,  $\nu<\frac {l-r}2$
for  $H=SO_0(l, r)$,
and $\nu <l-r +\frac 32$
for  $H=Sp(l, r)$.
If
$k\in \mathbb N$  such that
\begin{equation}
\label{cond-k-r}
0\le k <\frac 14(\frac{l-r} 2-\nu)
\end{equation}
for $H=SO_0(r, r+b)$
and
\begin{equation}\label{cond-k-h}
0\le k <\frac 18( 3+ 2(l-r) -2\nu)
\end{equation}
for $H=Sp(r, r+b)$, then
the spherical function $\phi_{\lam}$
with 
\begin{equation}
  \label{eq:lam-v-k}
i\lam=\begin{cases}
r(\nu +2k) -r( l-r +\frac 12(r-1)), &H=SO_0(l, r)\\
2r(\nu +2k) -r(3+l-r +8(r-1)), &H=Sp(l, r),
\end{cases}  
\end{equation}
is positive definite and appears
as a discrete component in the irreducible
representation of $(\mathcal H_{\nu}, G)$ under $H\subset G$.
\end{theo+}
\begin{proof}
The formula (\ref{eq:h-c-e}) for the Harish-Chandra
e-function implies that
the Poisson transform on $S$
is given by
$$
P_{\lam}f(z)= 
\int_S \frac
{
h(z, \bar z)^{\frac{\sig}2} 
}
{
h(z, \bar v) ^{\sig}
}
f(v) dv, \quad z\in X,
$$
which intertwines the induced representation $Ind_P^H(\lam)$ with
the regular action $X$. The functions
$h(z, z), h(z,  v)$ are polynomials in
 $z\in X$ and thus have holomorphic extension
to $z\in D$. Furthermore 
it is easy to see that
they have no zeros on $D$, by, e.g.,
the explicit formula of $h$ 
in terms of the determinant functions
\cite{Loos-bsd}.
Thus  $Pf(z)$
has a holomorphic extension on
$D$, still denoted by $P_{\lam}f(z)$. 
 For non negative
integer $k$ satisfying (\ref{cond-k-r})
and (\ref{cond-k-h})
we put $\sig: =\nu +2k$ with
the corresponding  $\lam$ as in Theorem 5.1;
this is explicitly computed  in (5.5).
We define, for $f\in L^2(S)$, 
\begin{equation}
\begin{split}
 T_{\nu, k}f(z):&= h(z, z)^{-\frac \nu 2}
P_{\lam}f(z) \\
&=h(z, z)^{k}
\int_S \frac{1}{h(z, v)^{\sig}} f(v) dv, \quad z\in D.
\end{split}
\end{equation}
 The transformation formula
(\ref{h-trans})
 of $h$  when restricted to $X$ is 
$$
h(gz,  \overline {gz})= 
J_{g}(z)^{\frac 1 p} 
h(z, \bar z)\overline{J_{g}(z)^{\frac 1 p} }
=J_{g}(z)^{\frac 2 p} 
h(z,  z), \quad g\in H, \quad z\in X,
$$
for the Jacobian $J_g$ is real-valued for $x\in X$. It's holomorphic extension
$h(z,  z)$  satisfies then
$$
h(gz,  {gz})= 
J_{g}(z)^{\frac 2 p} 
h(z, z).
$$
Thus $T_{\nu, k} $  intertwines
the induced representation 
 $Ind_P^H(\lam)$ with $(\mathcal O(D), H, \pi_{\nu})$.
We prove that $T_{\nu, k}$
maps the function $1$ into
$\mathcal H_{\nu}$ when $k$ satisfies
the stated condition. The rest is proved as in Theorem 4.5.
We rewrite
$$
T_{\nu, k}1(z) =
h(z, z)^{k}F(z), \quad F(z):=\int_S \frac{1}{h(z, v)^{\sig}} f(v) dv, \quad z\in D.
$$
and we shall prove that the function $F(z)$ is 
in $\mathcal H_{\nu}$, so is $T_{\nu}1(z) $ since
$h(z, z)^{k}$ is a polynomial of the coordinate functions
and each of them defines a bounded multiplication
operator on $\mathcal H_{\nu}$ \cite{az-mult}.

The function $F(z)$ apart from the factor of $h(z, z)^{\frac \sig 2}$
is the holomorphic extension of the spherical function
computed in Theorem 5.1. 
If $H=SO_0(l, r)$ we have, 
$$
F(z)=\sum_{\m}
\frac{
(\frac{\nu +2k}2)_{\m}
(\frac{\nu +2k +1}2)_{\m}
}
{
(\frac a2 (r-1) + b +\frac 12 )_{\m}
}
\frac
{\pi_{\m}
}
{(q)_{\m}}
p_{\n}(z).
$$
Its norm square in $\mathcal H_{\nu}$ is,
by Proposition 3.5, the hypergeometric function
\begin{equation}
\label{n-s-f}
\Vert F\Vert_{\mathcal H_{\nu} }^2=
 {}_4F_3{}(\frac{\sig}2,
\frac{\sig}2, \frac{\sig +1}2,
\frac{\sig +1}2; \frac{\nu}2,
\frac{\nu+1}2, \frac a2 (r-1) +b +\frac 12; 1^r).  
\end{equation}
The series is convergent if
$$
2\sig +1 -(\nu +1+b +\frac a2(r-1)) <-\frac a2 (r-1) 
$$
equivalently (recalling $\sig=\nu + 2k$)
$$
\nu +4k < \frac {l-r}2
$$
according to Lemma \ref{conv-le}, 
which is guaranteed by our assumption on $l, r$ and $\nu$. Thus
$F$ is in $\mathcal H_{\nu}$, so 
is $T_{\nu}1 $.

If  $H=Sp(l, r)$, which is of Type BC with $\n=(\m, \m)$, $a=4=2a^\prime$, the function $F$ is
$$
F(z)=\sum_{\n=(\m, \m)}\frac{(\sig)_{\m} (\sig-1)_{\m}  }
{(\frac{a}2(r-1)  +\frac{\iota +2b}2)  } p_{\n}(z),
$$
and
\begin{equation}
\label{n-s-f-2}
\Vert F\Vert_{\mathcal H_{\nu}}^2
={}_4 F_3(\sig, \sig, \sig-1, \sig-1; \frac a2(r-1)+b+\frac {\iota }2, \nu, \nu-1; 1^r),
\end{equation}
whose convergence is again determined by Lemma 4.1.
The condition on the convergence is $l-r>\nu + 4k -\frac 32 >
\frac{a^\prime} 2(r^\prime -1)-\frac 32 =
(2r -1)-\frac 32$, namely $l-r \ge 2r-2=2(r-1)$. The condition
on $k\ge 0$ is then obtained accordingly.
\end{proof}

\subsection{Discrete components of $(\mathcal H_\nu, G)$
under $H$ for $\nu=\frac {a^\prime} 2(j -1)$ being
a singular Wallach point}

Let $H=SO_0(l, r)$
or  $H=Sp(l, r)$
 be as in the previous subsection. 
We fix $\nu=\frac {a^\prime} 2(j -1)
=j-1$, $2\le j\le r^\prime$ be a singular Wallach point,
where $r^\prime=r$ and respectively
$r^\prime =2r$.
We consider the operator $f\mapsto T_{\nu}f:=T_{\nu, 0}f$
defined in (5.6), and the image $F=T_{\nu, 0}1$ of the constant
function $1$.
The norm square $\Vert F\Vert_{\mathcal H_{\nu}}^2$
is again a series as in (\ref{n-s-f})
and (\ref{n-s-f-2} ).  The condition for the convergence for the group
$SO_0(l, r)$ is $l-r>j-1$ while as for the
group $Sp(l, r)$  is $l-r >j-1-\frac 32$, namely
$l-r \ge j-2$. Thus we have the following

\begin{theo+} Let $H$ be 
group $SO_0(l, r)$
or $Sp(l, r)$
 and  let $\nu=j-1$, $2\le j\le r^\prime$, be
a singular Wallach point of $G$.
Suppose  $l-r >(j-1)$ for the group
 $H=SO_0(r, r+b)$, 
and  $l-r\ge j-2$ 
 for $H=Sp(r, r+b)$.
The spherical function $\phi_{\lam}$,
with  $\lam$ determined by $(\nu, k):=(j-1, 0)$ as
in (5.5), is positive definite and appears
as a discrete component in the irreducible
representation of $(\mathcal H_{\nu}, G)$ under $H\subset G$.
\end{theo+}

Note that in Theorems 4.5 and 5.3 we have only taken
the operator $T_{\nu, k}$ with $k=0$
when  $\nu$ is a singular Wallach point.
 It would be interesting to refine 
the definition of $T_{\nu, k}$ so that 
we get  also finitely many discrete components;
indeed
there are finitely many discrete components
in the branching of $\pi_{\nu}\otimes \overline{\pi_{\nu}}$
when $\nu=1$ is the last Wallach point [39].

\medskip

\bibliographystyle{amsplain}
\def\cprime{$'$} \newcommand{\noopsort}[1]{} \newcommand{\printfirst}[2]{#1}
  \newcommand{\singleletter}[1]{#1} \newcommand{\switchargs}[2]{#2#1}
  \def\cprime{$'$} \def\cprime{$'$} \def\cprime{$'$}
\providecommand{\bysame}{\leavevmode\hbox to3em{\hrulefill}\thinspace}
\providecommand{\MR}{\relax\ifhmode\unskip\space\fi MR }
\providecommand{\MRhref}[2]{%
  \href{http://www.ams.org/mathscinet-getitem?mr=#1}{#2}
}
\providecommand{\href}[2]{#2}

\end{document}